\documentclass[12pt]{amsart}
\usepackage{amssymb}
\usepackage[all]{xy}
\usepackage{tikz}

\usepackage{amscd}
\usepackage{amsmath}
\input xy
\xyoption{all}

\newtheorem{lemma}{Lemma}[section]

\newtheorem{theorem}[lemma]{Theorem}
\newtheorem{proposition}[lemma]{Proposition}
\newtheorem{construction}[lemma]{Construction}

\theoremstyle{definition}
\newtheorem{remark}[lemma]{Remark}
\newtheorem{notation}[lemma]{Notation}
\newtheorem{definition}[lemma]{Definition}

\DeclareMathOperator{\Mod}{Mod}
\DeclareMathOperator{\modd}{mod}

\DeclareMathOperator{\Hom}{Hom}

\DeclareMathOperator{\Ext}{Ext}
\DeclareMathOperator{\Lim}{Lim}

\DeclareMathOperator{\Ker}{Ker}
\DeclareMathOperator{\Coker}{Coker}
\DeclareMathOperator{\Cogen}{Cogen}
\DeclareMathOperator{\cogen}{cogen}
\DeclareMathOperator{\Gen}{Gen}
\DeclareMathOperator{\Imm}{Im}
\DeclareMathOperator{\brick}{brick}
\DeclareMathOperator{\Brick}{Brick}
\DeclareMathOperator{\sbrick}{sbrick}
\DeclareMathOperator{\Sbrick}{Sbrick}
\DeclareMathOperator{\Ann}{Ann}
\DeclareMathOperator{\D}{D}
\DeclareMathOperator{\F}{F}
\DeclareMathOperator{\ts}{t}
\DeclareMathOperator{\fs}{f}
\DeclareMathOperator{\torf}{torf}
\DeclareMathOperator{\Torf}{Torf}
\DeclareMathOperator{\HH}{H}
\DeclareMathOperator{\HRS}{HRS}
\DeclareMathOperator{\Prod}{Prod}
\DeclareMathOperator{\pprod}{prod}
\DeclareMathOperator{\Ind}{Ind}
\DeclareMathOperator{\Rad}{Rad}
\DeclareMathOperator{\Cosilt}{Cosilt}
\DeclareMathOperator{\tilt}{tilt}
\DeclareMathOperator{\proj}{proj}

\newtheorem*{theorem 0*}{Theorem}
\newtheorem*{theorem a*}{Theorem A}
\newtheorem*{theorem b*}{Theorem B}

\newcounter{diagram}
\numberwithin{diagram}{section}

\begin{document}
	
	\title{Semibrick-cosilting correspondence}
	
	\author{Ramin Ebrahimi}
	\address{School of Mathematics, Institute for Research in Fundamental Sciences (IPM), P.O. Box: 19395-5746, Tehran, Iran}
	\email{ramin.ebrahimi1369@gmail.com / r.ebrahimi@ipm.ir}
	
	\author{Alireza Nasr-Isfahani}
	\address{Department of Pure Mathematics\\
		Faculty of Mathematics and Statistics\\
		University of Isfahan\\
		P.O. Box: 81746-73441, Isfahan, Iran\\ and School of Mathematics, Institute for Research in Fundamental Sciences (IPM), P.O. Box: 19395-5746, Tehran, Iran}
	\email{nasr$_{-}$a@sci.ui.ac.ir / nasr@ipm.ir}

	\subjclass[2010]{{18E10}, {18E20}, {18E99}}
	
	\keywords{Brick, Semibrick, Cosilting module,  Support $\tau$-tilting module}

	\begin{abstract}
		Let $\Lambda$ be a finite dimensional algebra. In this paper we show that there is a natural bijection between cosilting modules in $\Mod \Lambda$ and semibricks in $\Mod \Lambda$ satisfying some condition. Also this bijection restricts to a bijection between all semibricks in $\modd \Lambda$ and a certain subclass of cosilting modules. These bijections are generalizations of Asai's correspondence \cite{A} between support $\tau^-$-tilting modules and right finite semibricks.
	\end{abstract}
	
	\maketitle


	\section{Introduction}
	Let $\Lambda$ be a finite dimensional algebra. Tilting theory of the category of finite dimensional $\Lambda$-modules has its origin in the reflection functors of Bernstein, Gelfand and Ponoromev \cite{BGP}. In the following years (after the seminal work of Auslander, Platzek and Reiten \cite{APR}) the abstract definition of a tilting module was given by Brenner and Butlere \cite{BB}. They introduced so called tilting functors that capture reflection functors of Bernstein, Gelfand and Ponoromev.
	Let $T=\bigoplus_{i=1}^nT_i\in \modd \Lambda$ be a tilting module, with $T_i$'s indecomposable and pairwise non-isomorphic. It is well known that we can replace $T_j$ with another indecomposable module $T_i'$ and obtain a new tilting module if and only if $\frac{T}{T_j}$ is a faithful, and if $\frac{T}{T_j}$ is a not faithful, we cannot replace $T_j$. This replacement, if it is possible, is called the {\em mutation} of $T$ at $T_j$.
	
	$\tau$-tilting theory was introduced in \cite{AIR}, in order to complete the class of tilting modules with respect to mutation. Indeed, the authors introduced the class of support $\tau$-tilting modules, containing tilting modules, such that the operation of mutation is possible in any indecomposable direct summand. One of the main theorems of \cite{AIR} is a natural bijection between support $\tau$-tilting modules (support $\tau^-$-tilting modules) and functorially finite torsion classes (torsion free classes) in $\modd \Lambda$. 
	
	Later Asai \cite{A} showed that for each support $\tau^-$-tilting module $M$, we can associate a set of bricks in $\modd\Lambda$, that he called it a semibrick, and he proved the following theorem.
	\begin{theorem 0*}
		Let $\Lambda$ be a finite dimensional algebra. Then there exists a bijection between support $\tau^-$-tilting modules and semibricks in $\modd \Lambda$ that generate a functorially finite torsion free class in $\modd\Lambda$.
	\end{theorem 0*}
	Our aim in this paper is to extend the above theorem in such a way that allows for the parametrization of all semibricks. In order to do so, we need infinitely generated analogous of support $\tau$-tilting modules.
	Infinitely generated analogous of support $\tau$-tilting modules (and their dual), silting (and cosilting) modules where introduced and studied by many authors, mainly by Lidia Angeleri H{\"u}gel and her co-authors \cite{AH,AMV,BP}.
	
	Our first main result in this paper is the following Theorem.
	
	\begin{theorem a*}(Theorem \ref{t3.4})
		Let $\Lambda$ be a finite dimensional algebra. Then there exists a bijection between equivalence classes of cosilting modules and bordered semibricks (see Definition \ref{d3.3}) in $\Mod \Lambda$. 
	\end{theorem a*}
	
	A torsionfree class $\fs\subseteq\modd \Lambda$ is called widely generated if there is a wide subcategory $\mathcal{W}\subseteq\modd \Lambda$ such that $\fs$ is the smallest torsion free class containing $\mathcal{W}$. Restricting the above theorem to semibricks in $\modd\Lambda$ we obtain the following result.
	
	\begin{theorem b*}(Theorem \ref{t3.10})
		Let $\Lambda$ be a finite dimensional algebra. Then there exists a bijection between equivalence classes of cosilting modules that cogenerate a widely generated torsion free class in $\modd \Lambda$ and semibricks in $\modd \Lambda$.
	\end{theorem b*}
	
	Support $\tau^-$-tilting modules are certainly cosilting modules, and so we can recover Asai's correspondence \cite{A} from the above result (see Theorem \ref{t3.13}).

	\subsection{Notation}
	Throughout this paper, $\Lambda$ is a finite dimensional algebra. We denote by $\Mod\Lambda$ ($\modd\Lambda$) the category of all (finitely presented) right $\Lambda$-modules, and by $\tau$, the Auslander-Reiten translation. Also we denote by $\D(\Lambda)$ the unbounded derived category of $\Mod\Lambda$.
	
	Let $\mathcal{A}$ be an additive category and $\mathcal{X}$ be a class of objects in $\mathcal{A}$. By $\Prod(\mathcal{X})$ ($\pprod(\mathcal{X})$) we mean the subcategory of all objects in $\mathcal{A}$ that are direct summand of a (finite) product of objects in $\mathcal{X}$. Moreover, if $\mathcal{A}$ is abelian, we denote by $\Cogen(\mathcal{X})$ ($\cogen(\mathcal{X})$) the subcategory of all objects that are subobject of an object in $\Prod(\mathcal{X})$ ($\pprod(\mathcal{X})$).
	
	Let $\mathcal{X}$ be a subcategory of $\D(\Lambda)$ and $I$ be a subset of integers. We set
	\begin{equation}
	\mathcal{X}^{\bot_I}:=\{M\in \D(\Lambda) \mid \Hom(X,M[i])=0, \forall i\in I \;\text{and}\; \forall X\in \mathcal{X}\}.\notag
	\end{equation}
	$^{\bot_I}\mathcal{X}$ is defined dually. As usual we denote the interval $\{i\in \mathbb{Z}\mid i\geq 0 \}$ by $\geq 0$ and similarly the interval $\{i\in \mathbb{Z}\mid i\leq 0 \}$ by $\leq 0$.
	
	In a similar way, for a subcategory $\mathcal{X}\subseteq \Mod\Lambda$ and a subset of non-negative integers $I$ we define $\mathcal{X}^{\bot_I}$ and $^{\bot_I}\mathcal{X}$. When $I$ consists of only one number, $I=\{n\}$, we will use the notation $\mathcal{X}^{\bot_n}$ and $^{\bot_n}\mathcal{X}$ instead of $\mathcal{X}^{\bot_{\{n\}}}$ and $^{\bot_{\{n\}}}\mathcal{X}$.
	
	Assume that $\mathcal{X}\subseteq \Mod\Lambda$ (resp., $\mathcal{X}\subseteq \modd\Lambda$). Then we denote by $\F(\mathcal{X})$ (resp., $\tilde{\F}(\mathcal{X})$) the smallest torsion free class in $\Mod\Lambda$ (resp., $\modd\Lambda$) which is contains $\mathcal{X}$.

	
	\section{preliminaries}
	In this section we recall the definitions of torsion pairs, cosilting modules, and their relation with each other. Also we recall some results about bricks, and the relation with torsion classes.
	
	\begin{definition}\label{d2.1}
		Let $\mathcal{A}$ be an abelian category. A pair of full subcategories $(\mathcal{T},\mathcal{F})$ is called a {\em torsion pair} if it satisfies the following conditions.
		\begin{itemize}
			\item[(1)] $\Hom_{\mathcal{A}}(T,F)=0$ for every objects $T\in \mathcal{T}$ and $F\in \mathcal{F}$.
			\item[(2)] For any object $A\in \mathcal{A}$, there is a short exact sequence 
			\begin{equation}
			0\rightarrow T_A\rightarrow A\rightarrow F_A\rightarrow 0 \notag
			\end{equation}
			with $T_A\in \mathcal{T}$ and $F_A\in \mathcal{F}$.
		\end{itemize}
		In this case $\mathcal{T}$ is called the {\em torsion class} and $\mathcal{F}$ is called the {\em torsion free class}.
	\end{definition}
	
	Torsion pairs in $\modd \Lambda$, ordered by inclusion of their torsion class, forms a complete lattice \cite{DIRRH}. This lattice had been topic of many researches \cite{AIR,DIJ,BKZ,E}. The characterization of functorially finite torsion pairs using $\tau$-tilting theory, is one of the first significant results in this direction. Let recall the definition of functorially finite subcategory.
	
	\begin{definition}\label{d2.2}
		Let $\mathcal{A}$ be an additive category and $\mathcal{B}$ be a full subcategory of $\mathcal{A}$. 
		\begin{itemize}
			\item[(1)] 	$\mathcal{B}$ is called {\em covariantly finite in $\mathcal{A}$} (or a {\em preenveloping subcategory} of $\mathcal{A}$), if for every $A\in \mathcal{A}$ there exist an object $B\in\mathcal{B}$ and a morphism $f : A\rightarrow B$ such that, for all $B'\in\mathcal{B}$, the sequence of abelian groups $\Hom_\mathcal{A}(B, B')\rightarrow \Hom_\mathcal{A}(A, B')\rightarrow 0$ is exact. Such a morphism $f$ is called a {\em left $\mathcal{B}$-approximation of $A$} (or a {\em $\mathcal{B}$-preenvelope of $A$}).
			\item[(2)] A left $\mathcal{B}$-approximation $f: A\rightarrow B$ is called a {\em minimal left $\mathcal{B}$-approximation} (or a {\em $\mathcal{B}$-envelope}), if any morphism $g:B\rightarrow B$ satisfying $gf=f$ is an isomorphism.
		\end{itemize}
	        The notions of {\em contravariantly
			finite subcategory of $\mathcal{A}$} (or{\em precovering subcategory}), {\em right $\mathcal{B}$-approximation} (or {\em $\mathcal{B}$-precover}) and {\em minimal right $\mathcal{B}$-approximation} (or{\em $\mathcal{B}$-cover}) are defined dually. A {\em functorially
			finite subcategory of $\mathcal{A}$} is a subcategory which is both covariantly and contravariantly finite
		    in $\mathcal{A}$.
	\end{definition}
	
	Let $(\mathcal{T},\mathcal{F})$ be a torsion pair in an abelian category $\mathcal{A}$ and $A\in \mathcal{A}$. By definition there is an exact sequence
	\begin{equation}
	0\rightarrow T_A\rightarrow A\rightarrow F_A\rightarrow 0 \notag
	\end{equation}
	with $T_A\in \mathcal{T}$ and $F_A\in \mathcal{F}$. From this exact sequence we can easily see that $T_A\rightarrow A$ is a $\mathcal{T}$-cover for $A$ and $A\rightarrow F_A$ is an $\mathcal{F}$-envelope for $A$. Thus $\mathcal{T}$ is a covering subcategory and $\mathcal{F}$ is an enveloping subcategory. For torsion pairs in $\modd\Lambda$ we have the following nice symmetry.
	
	\begin{theorem}$($\cite{S}$)$\label{t2.3}
		Let $(\ts,\fs)$ be a torsion pair in $\modd\Lambda$. Then $\ts$ is preenveloping if and only if $\fs$ is precovering.
	\end{theorem}
    A torsion pair $(\ts,\fs)$ in $\modd\Lambda$ is called functorially finite if $\ts$, or equivalently $\fs$, is functorially finite.
	
	Adachi, Iyama and Reiten in \cite{AIR} defined the notation of $\tau$-tilting modules and $\tau^-$-tilting modules in order to classify functorially finite torsion pairs.
	\begin{definition}\label{d2.4}
		Let $M\in \modd \Lambda$,
		\begin{itemize}
			\item[(1)] $M$ is called {\em $\tau$-rigid} (\em {$\tau^-$-rigid}) if $\Hom(M,\tau M)=0$ ($\Hom(\tau^- M,M)=0$).
			\item[(2)] $M$ is called a \em {$\tau$-tilting} (\em {$\tau^-$-tilting}) module if $M$ is $\tau^-$-rigid and $|M|=|\Lambda|$, where $|M|$ denotes the number of non-isomorphic indecomposable direct summands of $M$.
			\item[(3)] $M$ is called a \em {support $\tau$-tilting} (\em {support $\tau^-$-tilting}) module if $M$ is $\tau$-tilting ($\tau^-$-tilting) module over its support (i.e. there is an idempotent $e\in \Lambda$ such that $M$ is a $\tau$-tilting module in $\modd\Lambda/\langle e\rangle$).
		\end{itemize}
	\end{definition}
	
	For an arbitrary algebra $\Lambda$, support $\tau$-tilting $\Lambda$-modules are in bijection with many important objects, such as $2$-term silting complexes and bounded co-$t$-structures in $K^b(\proj \Lambda)$, simple minded collections and bounded $t$-structures in $\D^b(\Lambda)$ with a length heart. 
	
	One of the main results of \cite{AIR} is the following.
	\begin{theorem}\label{t2.5}
		There is a bijection between
		\begin{itemize}
			\item[(1)] Support $\tau$-tilting $\Lambda$-modules (Support $\tau^-$-tilting $\Lambda$-modules), and
			\item[(2)] Functorially finite torsion (torsion free) classes in $\modd \Lambda$.
		\end{itemize}
	\end{theorem}
	
	In this correspondence a support $\tau^-$-tilting $\Lambda$-module $M$ sends to $\cogen(M)$, and conversely a functorially finite torsion free class $\fs$ of $\modd \Lambda$ is send to the direct sum of all pairwise non-isomorphism $\Ext$-injective objects of $\fs$. So every functorially finite torsion free class as an exact category, induced from $\modd \Lambda$, is determined uniquely by it's injective objects. Nevertheless non-functorially finite torsion free classes have no injective object in general. So we cannot classify non-functorially finite torsion free classes using injectives (and $\tau^-$-tilting modules).
	There exists a method to address this issue.
	First let recall the following theorem due to Crowly Bowey \cite{CB}. Recall that a subcategory of $\Mod\Lambda$ is called {\em definable}, if it is closed under pure subobjects, products and direct limits. Because any torsion free class in $\Mod\Lambda$ is closed under subobjects and products, a torsion free class $\mathcal{F}$ in $\Mod\Lambda$ is definable if and only if it is closed under direct limits (i.e. $\mathcal{F}=\underrightarrow{\Lim}\mathcal{F}$).
	
	\begin{theorem} $($\cite{CB}$)$\label{t2.6}
		There is a bijection between
		\begin{itemize}
			\item[(1)] The set of all torsion pairs in $\modd \Lambda$.
			\item[(2)] The set of all torsion pairs in $\Mod \Lambda$, with definable torsion free class.
		\end{itemize}
		This bijection associates to a torsion pair $(\ts,\fs)$, its limit closure $(\underrightarrow{\Lim}\ts,\underrightarrow{\Lim}\fs)=(\Gen\ts, \ts^{\bot})$.
		The inverse map sends a torsion pair $(\mathcal{T},\mathcal{F})$ with $\mathcal{F}=\underrightarrow{\Lim}\mathcal{F}$, to the restriction $(\mathcal{T}\cap \modd \Lambda,\mathcal{F}\cap \modd \Lambda)$.
	\end{theorem}
	
	So by the above theorem, in order to classify all torsion pairs in $\modd \Lambda$ it is enough to classify all torsion pairs in $\Mod\Lambda$ with definable torsion free class. And as we will see in Theorem \ref{t2.9} this torsion pairs are in bijection with large analogous of support $\tau^-$-tilting modules.
	
	Cosilting modules were defined in \cite{BP} as the dual notation of {\em silting modules} \cite{AMV}.
	Let $\xi:Q_0\rightarrow Q_1 $ be a homomorphism between injective $\Lambda$-modules. Define
	\begin{equation*}
	\beta_{\xi}=\{X\in \Mod \Lambda |\Hom_{\Lambda}(X,\xi) \;\text{is an epimorphism}\}.
	\end{equation*}
	\begin{definition}\label{d2.7}
		$C\in \Mod \Lambda$ is called
		\begin{itemize}
			\item[(I)] \em {partial cosilting module} (with respect to $\xi$), if there exists an injective copresentation
			\begin{equation}\label{IC}
				0\rightarrow C\overset{f}{\rightarrow}Q_0\overset{\xi}{\rightarrow} Q_1
			\end{equation}
			of $C$ such that
			\begin{itemize}
				\item[(a)] $C\in \beta_{\xi}$, and
				\item[(b)] $\beta_{\xi}$ is closed under products, or equivalently $\beta_{\xi}$ is a torsion free class.
			\end{itemize}
			\item[(II)] \em {cosilting module} (with respect to $\xi$), if there exists an injective copresentation \eqref{IC} of $C$ such that $\Cogen(C)=\beta_{\xi}$.
		\end{itemize}
		Two cosilting modules $C$ and $C'$ are said to be equivalent if $\Prod(C)=\Prod(C')$ or equivalently $\Cogen(C)=\Cogen(C')$. We denote by $\Cosilt\Lambda$, the set of equivalence classes of cosilting $\Lambda$-modules.
	\end{definition}
	
	Recall that a cotilting module (of projective dimension 1) in $\Mod \Lambda$ is a module $C$ such that $\Cogen(C)=^{\bot_1}C$. The relation between silting modules and cotilting modules is given in the following proposition.
	
	\begin{proposition}\label{p2.8}
		\begin{itemize}
			\item[(1)] $($\cite{BP}$)$ A module $C$ is cotilting if and only if it is a cosilting module of injective dimension $\leq 1$.
			\item[(2)] $($\cite{ZW}$)$ The following conditions are equivalent for $C\in \Mod \Lambda$.
			\begin{itemize}
				\item[(a)] $C$ is a cosilting module.
				\item[(b)] $\Cogen(C)$ is a torsion free class and $C$ is a cotilting module over $\bar{\Lambda}=\frac{\Lambda}{\Ann{(C)}}$.
			\end{itemize}
		\end{itemize}
	\end{proposition}

	\begin{theorem}$($\cite{ZW}$)$\label{t2.9}
		There is a bijection between the equivalence classes of cosilting $\Lambda$-modules and torsion pairs with definable torsion free class in $\Mod\Lambda$. This bijection sends a cosilting module $C$ to the torsion pair $(^{\bot_0}C,\Cogen(C))$.
	\end{theorem}

   \begin{remark}\label{r2.10}
   	Let $\torf\Lambda$ be the set of torsion free classes in $\modd\Lambda$ and $\Torf_d$ be the set of definable torsion free classes in $\Mod\Lambda$.
   	By Theorem 2.6 and Theorem 2.9 we have the following bijections:
   	\begin{align*}
   		&\Cosilt\Lambda &\longleftrightarrow &\Torf_d\Lambda \longleftrightarrow &\torf\Lambda\\
   		&C&\longmapsto &\Cogen(C)\longmapsto &\Cogen(C)\cap\modd\Lambda
   	\end{align*}
   \end{remark}
	
	The lattice theory of torsion pairs in $\modd\Lambda$ is controlled by specific modules, referred to as bricks.
	
	\begin{definition}\label{d2.11}
		let $\mathcal{A}$ be an abelian category.
		\begin{itemize}
			\item[(1)] An object $S\in \mathcal{A}$ is called a {\em brick}, if the endomorphism ring of $S$ is a division ring (i.e. every nonzero map $f:S\rightarrow S$ is an isomorphism).
			\item[(1)] A set $\mathcal{S}$ of bricks is called a {\em semibrick} if $\Hom(S_1,S_2)=0$ for any $S_1\neq S_2\in \mathcal{S}$.
		\end{itemize}
	\end{definition}
	
	\begin{notation}\label{n2.12}
		Let $\Lambda$ be a finite dimensional algebra.
		\begin{itemize}
			\item[(1)] We denote by $\brick \Lambda$ (respectively $\Brick \Lambda$) the set of all bricks in $\modd \Lambda$ (respectively in $\Mod \Lambda$).
			\item[(2)] We denote by $\sbrick \Lambda$ (respectively $\Sbrick \Lambda$) the set of all semibricks in $\modd \Lambda$ (respectively in $\Mod \Lambda$).
		\end{itemize}
	\end{notation}
	
	Bricks and semibricks are closely related to the concept of torsion free, almost torsion modules \cite{AHL}.
	
	\begin{definition}$($\cite[Definition 3.1]{AHL}$)$\label{d2.13}
		let $(\mathcal{T},\mathcal{F})$ be a torsion pair in abelian category $\mathcal{A}$. $S\in \mathcal{A}$ is called {\em torsion free, almost torsion} if the following conditions are satisfied.
		\begin{itemize}
			\item[(1)] $S\in \mathcal{F}$.
			\item[(2)] Every proper quotient of $S$ is contained in $\mathcal{T}$.
			\item[(3)] For every $F\in \mathcal{F}$ and every short exact sequence $0\rightarrow S\rightarrow F\rightarrow C\rightarrow 0$, $C\in \mathcal{F}$.
		\end{itemize}
	\end{definition}
	
	We mention that torsion free, almost torsion objects for torsion pairs in $\modd\Lambda$ were studied by Barnard, Carroll, and Zhu under the name {\em minimal extending modules} \cite{BKZ}.
	 
	The following proposition gives an equivalent definition for torsion free, almost torsion modules.
	\begin{proposition}$($\cite[\S 3]{S}$)$\label{p2.14}
		let $(\mathcal{T},\mathcal{F})$ be a torsion pair in abelian category $\mathcal{A}$ and $S\in \mathcal{A}$.
		\begin{itemize}
			\item[(1)] $S$ satisfies the condition $(2)$ of Definition \ref{d2.13} if and only if any nonzero map $S\rightarrow F$ with $F\in \mathcal{F}$ is monomorphism.
			\item[(2)] Let $S$ satisfies the conditions $(1)$ and $(2)$ of  Definition \ref{d2.13}. Then $S$ satisfies the condition $(3)$ of Definition \ref{d2.13} if and only if for every non-split short exact sequence $0\rightarrow S\rightarrow E\rightarrow T\rightarrow 0$ with $T\in \mathcal{T}$, $E\in \mathcal{T}$.
		\end{itemize}
	\end{proposition}
	
	In the following easy proposition we state the relation between torsion free, almost torsion objects and bricks.
	
	\begin{proposition}\label{p2.15}
		let $(\mathcal{T},\mathcal{F})$ be a torsion pair in abelian category $\mathcal{A}$.
		\begin{itemize}
			\item[(1)] Every torsion free, almost torsion object is a brick.
			\item[(2)] The set of isomorphism classes of all torsion free, almost torsion objects forms a semibrick.
		\end{itemize}
		\begin{proof}
			Let $S_1,S_2\in \mathcal{F}$ be two torsion free, almost torsion objects and $f:S_1\rightarrow S_2$ be a nonzero morphism. By Proposition \ref{p2.14} $f$ is a monomorphism. Also by definition $\Coker(f)\in \mathcal{F} \cap \mathcal{T}=0$. Thus $f$ is an isomorphism.
		\end{proof}
	\end{proposition}
	
	\subsection{Torsion free, almost torsion modules in a cotilting torsion pair}
	Let $C$ be a cotilting module. Denote by $t_C=(^{\bot_0}C,\Cogen(C))$, the torsion pair cogenerated by $C$. In this subsection, following \cite{AHL}, we give an explicit description of torsion free, almost torsion modules. To this end we need $t$-structures \cite{BBD}.
	
	\begin{definition}\label{d2.16}
		A {\em $t$-structure} in $\D(\Lambda)$ is a pair of subcategories $(\mathcal{X},\mathcal{Y})$, satisfying the following conditions.
		\begin{itemize}
			\item[(1)] $\Hom(\mathcal{X},\mathcal{Y})=0$.
			\item[(2)] For any object $M\in \D(\Lambda)$, there exist a triangle $X_M\rightarrow M\rightarrow Y_M\rightarrow X_M[1]$, where $X_M\in \mathcal{X}$ and $Y_M\in \mathcal{Y}$.
			\item[(3)] $\mathcal{X}[1]\subseteq \mathcal{X}$.
		\end{itemize}
	\end{definition}
	The heart of $t$-structure $(\mathcal{X},\mathcal{Y})$ is the subcategory $\mathcal{H}=\mathcal{X}[-1]\cap\mathcal{Y}$. It was shown in \cite{BBD} that the heart of any $t$-structure is an abelian category.
	
	Happel, Reiten and Smalø showed in \cite{HRS} that any torsion pair $t=(\mathcal{T},\mathcal{F})$ in $\Mod \Lambda$ induces a $t$-structure $(\mathcal{X}_t,\mathcal{Y}_t)$ in $\D(\Lambda)$, where
	
	\begin{align}
	&\mathcal{X}_t=\{X\in \D(\Lambda)\mid\HH^0(X)\in \mathcal{T} \space\text{and}\space \HH^i(X)=0, \forall i\geq 1\} \notag \\
	&\mathcal{Y}_t=\{X\in \D(\Lambda)\mid\HH^0(X)\in \mathcal{F} \space\text{and}\space \HH^i(X)=0, \forall i\leq -1\}\notag
	\end{align}
	And the heart of this $t$-structure, denoted by $\mathcal{H}_t$, admits a torsion pair $(\mathcal{F},\mathcal{T}[-1])$. $(\mathcal{X}_t,\mathcal{Y}_t)$ is called $\HRS$-tilt of $t$.
	
	In the following proposition we collect some of the basic properties of $\HRS$-tilt of a cotilting torsion pair.
	
	\begin{proposition}$($\cite{CGF}$)$\label{p2.17}
		Let $C$ be a cotilting module and $t=t_C=(^{\bot_0}C,\Cogen(C))$ be
		the torsion pair cogenerated by $C$.
		\begin{itemize}
			\item[(1)] $\mathcal{H}_t$ is a Grothendieck category.
			\item[(2)] $C$ is an injective cogenerator for $\mathcal{H}_t$.
		\end{itemize}
	\end{proposition}

Since $(\Cogen(C),^{\bot_0}C[-1])$ is a torsion pair, simple objects of $\mathcal{H}_t$ are of two forms. Either they belong to $\Cogen(C)$ or they belong to $^{\bot_0}C[-1]$. Following \cite{AHL} we call simple objects in $\Cogen(C)$, the torsion simples. The following theorem characterizes injective envelopes of torsion simple objects. Recall that a {\em left almost split morphism} $f:M\rightarrow N$ in an additive category is a morphism that is not split monomorphism, and any other morphism $f':M\rightarrow N'$ which is not split monomorphism factor trough $f$. When such a factorization is unique, $f$ is called {\em strong left almost split}.

\begin{theorem}$($\cite{AHL}$)$\label{t2.18}
	Let $t=t_C=(^{\bot_0}C,\Cogen(C))$ be the torsion pair cogenerated by a cotilting module $C$. The torsion simple objects in $\mathcal{H}_t$ are exactly torsion free almost torsion modules and the following statements are equivalent for a module $X$.
	\begin{itemize}
		\item[(1)] $X$ is isomorphic to the injective envelope of a torsion simple $S\in \mathcal{H}_t$.
		\item[(2)] There exists a short exact sequence
		\begin{equation}
		0\rightarrow S\overset{a}{\rightarrow} X\overset{b}{\rightarrow} \bar{X}\rightarrow 0 \notag,
		\end{equation}
		in $\Mod\Lambda$, where $S$ is torsion free, almost torsion, $b$ is a strong left almost split morphism in $\Cogen(C)$ and $a$ is a $(\Cogen C)^{\bot_1}$-envelope of $S$.
	\end{itemize}
\end{theorem}

	\begin{definition}\label{d2.19}
		let $C\in \Mod\Lambda$ be a cosilting module. We define $\Ind_C$ to be the set of equivalence classes of indecomposable objects in $\Prod(C)$.
	\end{definition}
	
	The following lemma is essential in the sequel.
	\begin{lemma}$($\cite{ALS}$)$\label{l2.20}
		let $C\in \Mod\Lambda$ be a cosilting module. Put $C':=\prod_{X\in \Ind_C} X$. Then $C'$ is a cosilting module equivalent to $C$.
	\end{lemma}
	
	\begin{construction}\label{c2.21}
		Let $C\in \Mod\Lambda$ be a cosilting module. For any $X\in\Ind_C$ consider the canonical morphism
		\begin{equation}\label{phi}
		\phi_X:X\longrightarrow \prod_{Y\in\Ind_C}Y^{\Rad(X,Y)},
		\end{equation}
		where $\Rad(X,Y)=\{h\in \Hom_{\Lambda}(X, Y)|1_X-gh \space\text{ is invertible for any }\space  g\in\Hom_\Lambda(Y, X)\}$. Set $S_X:=\Ker(\phi_X)$ and $\mathcal{S}_C:=\{S_X\mid X\in\Ind_C\; \text{and}\; S_X\neq 0\}$. Also for any $X\in\Ind_C$ we have the following short exact sequence
		\begin{equation}\label{cses}
		0\rightarrow S_X\overset{a}{\rightarrow} X\overset{b}{\rightarrow} \bar{X}\rightarrow 0,
		\end{equation}
		where $\bar{X}=\Imm(\phi_X)$.
	\end{construction}
	
	We keep the notation of Construction \ref{c2.21} through the paper.
	
	\begin{proposition}\label{p2.22}
		Let $S_X\neq 0$. Then in the short exact sequence \eqref{cses} $b$ is a strong left almost split morphism in $\Cogen(C)$ and $a$ is a $(\Cogen C)^{\bot_1}$-envelope of $S_X$.
		\begin{proof}
			First we show that $b$ is a strong left almost split morphism in $\Cogen C$. Let $f:X\rightarrow F$ be a non-split monomorphism in $\Cogen(C)$. Using Lemma \ref{l2.20} we obtain a monomorphism $i:F\rightarrow \prod_{Y\in \Ind_C}Y^{I_Y}$ for some index sets $I_Y$. So we have the solid part of the following diagram.
			\begin{center}
				\begin{tikzpicture}
				\node (X0) at (-4,1) {$0$};
				\node (X1) at (-2,1) {$S_X$};
				\node (X2) at (0,1) {$X$};
				\node (X3) at (2,1) {$\bar{X}$};
				\node (X4) at (4,1) {$0$};
				\node (X5) at (0,-1) {$F$};
				\node (X6) at (0,-3) {$\prod_{Y\in \Ind_C}Y^{I_Y}$};
				\draw [->,thick] (X0) -- (X1) node [midway,left] {};
				\draw [->,thick] (X1) -- (X2) node [midway,above] {$a$};
				\draw [->,thick] (X2) -- (X3) node [midway,above] {$b$};
				\draw [->,thick] (X3) -- (X4) node [midway,left] {};
				\draw [->,thick] (X2) -- (X5) node [midway,left] {$f$};
				\draw [>->,thick] (X5) -- (X6) node [midway,left] {$i$};
				\draw [->,thick,dashed] (X3) -- (X5) node [midway,left] {$\bar{f}$};
				\end{tikzpicture}
			\end{center}
			By the construction of $S_X$, and the assumption that $f$ is not split monomorphism, $ifa=0$. Since $i$ is a monomorphism, $fa=0$. Then by the universal property of cokernel, there is a unique morphism $\bar{f}:\bar{X}\rightarrow F$ making the diagram commutative.
			
			Now we prove that $a$ is a $(\Cogen C)^{\bot_1}$-envelope of $S_X$. First we note that $X\in (\Cogen C)^{\bot_1}$. Let $g:S_X\rightarrow M$ be a morphism with $M\in (\Cogen C)^{\bot_1}$. By taking push out along $g$ we obtain the following commutative diagram with exact rows.
			\begin{center}
				\begin{tikzpicture}
				\node (X0) at (-4,1) {$0$};
				\node (X1) at (-2,1) {$S_X$};
				\node (X2) at (0,1) {$X$};
				\node (X3) at (2,1) {$\bar{X}$};
				\node (X4) at (4,1) {$0$};
				\node (X5) at (-4,-1) {$0$};
				\node (X6) at (-2,-1) {$M$};
				\node (X7) at (0,-1) {$E$};
				\node (X8) at (2,-1) {$\bar{X}$};
				\node (X9) at (4,-1) {$0$};
				\draw [->,thick] (X0) -- (X1) node [midway,left] {};
				\draw [->,thick] (X1) -- (X2) node [midway,above] {$a$};
				\draw [->,thick] (X2) -- (X3) node [midway,above] {$b$};
				\draw [->,thick] (X3) -- (X4) node [midway,left] {};
				\draw [->,thick] (X5) -- (X6) node [midway,left] {};
				\draw [->,thick] (X6) -- (X7) node [midway,left] {};
				\draw [->,thick] (X7) -- (X8) node [midway,left] {};
				\draw [->,thick] (X8) -- (X9) node [midway,left] {};
				\draw [->,thick] (X1) -- (X6) node [midway,left] {$g$};
				\draw [->,thick] (X2) -- (X7) node [midway,left] {};
				\draw [double,-,thick] (X3) -- (X8) node [midway,left] {};
				\end{tikzpicture}
			\end{center}
			By assumption $\Ext^1(\bar{X},M)=0$, so the bottom row splits. This implies that there is a morphism $t:X\rightarrow M$ such that $ta=g$. Thus $a$ is a $(\Cogen C)^{\bot_1}$-preenvelope for $S_X$. It remains to show that $a$ is left minimal. Let $f:X\rightarrow X$ be a morphism such that $fa=a$. If $f$ is not isomorphism, by the construction of $a:S_X\rightarrow X$, we have that $fa=0$. Thus $a=0$ that contradicts the assumption that $S_X\neq 0$.
		\end{proof}
	\end{proposition}
	
	\begin{proposition}\label{p2.23}
		Let $S_X\neq 0$. Then $S_X$ is a torsion free almost torsion module for the torsion pair $t_C=(^{\bot_0}C,\Cogen C)$.
		\begin{proof}
			Let $f:S_X\rightarrow F$ be a nonzero morphism with $F\in \Cogen (C)$. We must show that $f$ is a monomorphism. As in the proof of Proposition \ref{p2.22} there is a monomorphism $i:F\rightarrow \prod_{Y\in \Ind_C}Y^{I_Y}$. Since $\Ext^1(\bar{X},\prod_{Y\in \Ind_C}Y^{I_Y})=0$, there is a map $g:X\rightarrow \prod_{Y\in \Ind_C}Y^{I_Y}$ making the following diagram commutative.
			\begin{center}
				\begin{tikzpicture}
				\node (X0) at (-4,1) {$0$};
				\node (X1) at (-2,1) {$S_X$};
				\node (X2) at (0,1) {$X$};
				\node (X3) at (2,1) {$\bar{X}$};
				\node (X4) at (4,1) {$0$};
				\node (X5) at (-2,-1) {$F$};
				\node (X6) at (-2,-3) {$\prod_{Y\in \Ind_C}Y^{I_Y}$};
				\draw [->,thick] (X0) -- (X1) node [midway,left] {};
				\draw [->,thick] (X1) -- (X2) node [midway,above] {$a$};
				\draw [->,thick] (X2) -- (X3) node [midway,above] {$b$};
				\draw [->,thick] (X3) -- (X4) node [midway,left] {};
				\draw [->,thick] (X1) -- (X5) node [midway,left] {$f$};
				\draw [>->,thick] (X5) -- (X6) node [midway,left] {$i$};
				\draw [->,thick,dashed] (X2) -- (X6) node [midway,left] {$g$};
				\end{tikzpicture}
			\end{center}
			If $g$ is split monomorphism, then $ga$ is monomorphism and so $f$ is monomorphism. If $g$ is not split monomorphism, by Proposition \ref{p2.22} there is a map $\tilde{g}:\bar{X}\rightarrow \prod_{Y\in \Ind_C}Y^{I_Y}$ such that $g=\tilde{g}b$. Then $f=0$ that is a contradiction.
			
			Now we show that for any monomorphism $f:S_X\rightarrow F$ with $F\in \Cogen (C)$, we have that $\Coker(f)\in \Cogen(C)$. Put $G:=\Coker(f)$. Taking push out, we obtain the following commutative diagram with exact rows and columns.
			\begin{center}
				\begin{tikzpicture}
				\node (Y) at (-2,3) {$0$};
				\node (Z) at (0,3) {$0$};
				\node (X0) at (-4,1) {$0$};
				\node (X1) at (-2,1) {$S_X$};
				\node (X2) at (0,1) {$X$};
				\node (X3) at (2,1) {$\bar{X}$};
				\node (X4) at (4,1) {$0$};
				\node (X5) at (-4,-1) {$0$};
				\node (X6) at (-2,-1) {$M$};
				\node (X7) at (0,-1) {$E$};
				\node (X8) at (2,-1) {$\bar{X}$};
				\node (X9) at (4,-1) {$0$};
				\node (X10) at (-2,-3) {$G$};
				\node (X11) at (0,-3) {$G$};
				\node (X12) at (-2,-5) {$0$};
				\node (X13) at (0,-5) {$0$};
				\draw [->,thick] (Y) -- (X1) node [midway,left] {};
				\draw [->,thick] (Z) -- (X2) node [midway,left] {};
				\draw [->,thick] (X0) -- (X1) node [midway,left] {};
				\draw [->,thick] (X1) -- (X2) node [midway,above] {$a$};
				\draw [->,thick] (X2) -- (X3) node [midway,above] {$b$};
				\draw [->,thick] (X3) -- (X4) node [midway,left] {};
				\draw [->,thick] (X5) -- (X6) node [midway,left] {};
				\draw [->,thick] (X6) -- (X7) node [midway,left] {};
				\draw [->,thick] (X7) -- (X8) node [midway,left] {};
				\draw [->,thick] (X8) -- (X9) node [midway,left] {};
				\draw [->,thick] (X1) -- (X6) node [midway,left] {$f$};
				\draw [->,thick] (X2) -- (X7) node [midway,left] {$g$};
				\draw [double,-,thick] (X3) -- (X8) node [midway,left] {};
				\draw [->,thick] (X6) -- (X10) node [midway,left] {};
				\draw [->,thick] (X7) -- (X11) node [midway,left] {};
				\draw [->,thick] (X10) -- (X12) node [midway,left] {};
				\draw [->,thick] (X11) -- (X13) node [midway,left] {};
				\draw [double,-,thick] (X10) -- (X11) node [midway,left] {};
				\end{tikzpicture}
			\end{center}
			Because $\Cogen(C)$ is extension closed, $E\in \Cogen(C)$.
			We claim that $g$ is a split monomorphism. If not, by Proposition \ref{p2.22}, $g$ factor through $b$ and then $a=0$ which contradicts the assumption that $S_X\neq 0$. Thus $G$ is a direct summand of $E$ and so $G\in \Cogen(C)$.
		\end{proof}
	\end{proposition}
	
	\begin{theorem}\label{t2.24}
		$\mathcal{S}_C$ is the set of all torsion free, almost torsions for the torsion pair $t_C=(^{\bot_0}C,\Cogen C)$.
		\begin{proof}
			First we note that torsion free, almost torsion modules for $t_C=(^{\bot_0}C,\Cogen(C))$ as a torsion pair in $\Mod\Lambda$ are the same with torsion free, almost torsion modules for $t_C=(^{\bot_0}C,\Cogen(C))$ as a torsion pair in $\Mod\frac{\Lambda}{\Ann(C)}$. Thus by Proposition \ref{p2.8} without lose of generality we may assume that $C$ is a cotilting module.
			
			Let $S$ be a torsion free, almost torsion module for $t_C$. We must show that there is an indecomposable $X\in \Prod(C)$ such that $S= S_X$.
			By Theorem \ref{t2.18}, $S$ is a simple object in the cotilting heart $\mathcal{H}_t$. Thus there is an exact sequence
			\begin{equation}
			0\rightarrow S\overset{a}{\rightarrow} X\overset{b}{\rightarrow} \bar{X}\rightarrow 0 \notag,
			\end{equation}
		    where $a$ is a $(\Cogen(C))^{\bot_1}$-envelope and $b$ is a strong left almost split morphism in $\Cogen(C)$. We claim that $S=S_X$. Since in this case $S$ and $S_X$ are simple objects in the $\mathcal{H}_t$ and they have the same injective envelope $X$ by Theorem \ref{t2.18}, it is enough to show that $S_X\neq 0$. 
		    
		    If $S_X=0$ then $\phi_X:X\rightarrow \prod_{Y\in \Ind_C}Y^{\Rad(X,Y)}$ is a monomorphism. In the following commutative diagram
		    \begin{center}
		    	\begin{tikzpicture}
		    	\node (X0) at (-4,1) {$0$};
		    	\node (X1) at (-2,1) {$S$};
		    	\node (X2) at (0,1) {$X$};
		    	\node (X3) at (2,1) {$\bar{X}$};
		    	\node (X4) at (4,1) {$0$};
		    	\node (X5) at (0,-1) {$\prod_{Y\in \Ind_C}Y^{I_Y}$};
		    	\draw [->,thick] (X0) -- (X1) node [midway,left] {};
		    	\draw [->,thick] (X1) -- (X2) node [midway,above] {$a$};
		    	\draw [->,thick] (X2) -- (X3) node [midway,above] {$b$};
		    	\draw [->,thick] (X3) -- (X4) node [midway,left] {};
		    	\draw [->,thick] (X2) -- (X5) node [midway,left] {$\phi_X$};
		    	\draw [->,thick,dashed] (X3) -- (X5) node [midway,left] {$\bar{f}$};
		    	\end{tikzpicture}
		    \end{center}
	    because $b$ is a left almost split sequence in $\Cogen(C)$, the dashed arrow $\bar{f}$ exists that makes the diagram commutative. Then we have that $\phi_X a=0$ and since $\phi_X$ is a monomorphism, $a=0$. But this means that $S=0$ which is a contradiction.
		\end{proof}
	\end{theorem}
	
	\section{The main theorems}
	In this section we prove the main results of the paper. By Theorem \ref{t2.24} we have a well defined map.
	\begin{align*}
	\alpha:\Cosilt \Lambda&\longrightarrow \Sbrick \Lambda\\
	C&\longmapsto \mathcal{S}_C
	\end{align*}
	The following proposition tells that this map is injective.
	
	\begin{proposition}\label{p3.1}
		Two cosilting modules $C$ and $C'$ are equivalent if and only if $\mathcal{S}_C=\mathcal{S}_{C'}$.
		\begin{proof}
			Two cosilting modules $C$ and $C'$ are equivalent if and only if $\Cogen(C)=\Cogen(C')$. By \cite[Proposition 3.20]{S} we know that two definable torsion free class, are equal if and only if they have the same torsion free, almost torsion modules. Thus the result follow from Theorem \ref{t2.24}.
		\end{proof}
	\end{proposition}
	
	In order to get a bijection, we characterize semibricks that are in the image of $\alpha$.
	
	\begin{proposition}\label{p3.2}
		Let $C$ be a cosilting module. Then $\Cogen(C)=\underrightarrow{\Lim} \F(\mathcal{S}_C)$.
		\begin{proof}
			Note that $\Cogen(C)\cap \modd\Lambda\subseteq \F(\mathcal{S}_C)$, because any $M\in \Cogen(C)\cap \modd\Lambda$ is a finite filtration of simple objects (in the sense of \cite{E}) in the torsion free class $\Cogen(C)\cap \modd\Lambda$, and by \cite[Proposition 3.20]{S} any simple object in $\Cogen(C)\cap \modd\Lambda$ is isomorphic to a submodule of some torsion free, almost torsion module in $\mathcal{S}_C$. Thus
			\begin{center}
				$\Cogen(C)\cap \modd\Lambda\subseteq \F(\mathcal{S}_C)\subseteq \Cogen(C).$
			\end{center}
		    By taking limit closure we have that $\Cogen(C)=\underrightarrow{\Lim} \F(\mathcal{S}_C)$ by Theorem \ref{t2.9}.
		\end{proof}
	\end{proposition}
	
	Motivated by the above proposition we give the following definition.
	\begin{definition}\label{d3.3}
		We say that a semibrick $\mathcal{S}\in \Sbrick \Lambda$ is {\em bordered} if $\underrightarrow{\Lim} \F(\mathcal{S})$ is a definable torsion free class, with $\mathcal{S}$ as the set of torsion free, almost torsion modules. We denote by $\Sbrick_b \Lambda$ the set of all bordered semibricks in $\Mod\Lambda$.
	\end{definition}
	
	Now we can prove Theorem A.
	
	\begin{theorem}\label{t3.4}
		The map $\alpha:\Cosilt \Lambda\longrightarrow \Sbrick \Lambda$ gives a bijection
		\begin{align*}
		\alpha_1:\Cosilt \Lambda&\longrightarrow \Sbrick_b \Lambda\\
		C&\longmapsto \mathcal{S}_C
		\end{align*}
		between equivalence classes of cosilting modules and the set of bordered semibricks.
		\begin{proof}
			For any cosilting module $C$, by Proposition \ref{p3.2}, $\mathcal{S}_C$ is a bordered semibrick. Thus $\alpha$ defines a well defined map $\alpha_1:\Cosilt \Lambda\longrightarrow \Sbrick_b \Lambda$.
			
			Now let $\mathcal{S}$ be a bordered semibrick. Then $\underrightarrow{\Lim} \F(\mathcal{S})$ is a definable torsion free class. Thus by Theorem \ref{t2.9} there exists a cosilting module $C$ such that $\underrightarrow{\Lim} \F(\mathcal{S})=\Cogen(C)$. By Proposition \ref{p2.24}, $\mathcal{S}=\mathcal{S}_C=\alpha(C)$.
		\end{proof}
	\end{theorem}
	
	Now we go on to prove Theorem B.
	\begin{notation}\label{n3.5}
		Let $C\in \Mod\Lambda$ be a cosilting module. We denote by $(\ts_C,\fs_C)$ the associated torsion pair in $\modd \Lambda$ (i.e. $(\ts_C,\fs_C)=(^{\bot_0}C\cap\modd \Lambda, \Cogen(C)\cap \modd \Lambda)$).
	\end{notation}
	Recall that a {\em wide subcategory} of an abelian category is a subcategory closed under kernels, cokernels and extensions.
	\begin{definition}$($\cite{A}$)$\label{d3.6}
		A torsion free class $\fs\subseteq \modd\Lambda$ is called {\em widely generated} if there is a wide subcategory $\mathcal{W}\subseteq \modd\Lambda$ such that $\fs=\tilde{\F}(\mathcal{W})$.
		We denote by $\Cosilt_w\Lambda$ the set of equivalence classes of cosilting module $C$ such that the associated torsion free class $\fs_C=\Cogen(C)\cap \modd \Lambda$ is widely generated.
	\end{definition}

Let $\fs\subseteq \modd\Lambda$ be a torsion free class. Following \cite{IT,MS} we define
\begin{align*}
	\beta(\fs):=\{X\in \fs\mid \forall(h:X\rightarrow Y)\in \fs, \Coker(h)\in \fs\}.
\end{align*}

\begin{proposition}\label{p3.7}
	Let $(\ts,\fs)$ be a torsion pair in $\modd\Lambda$.
	\begin{itemize}
		\item[(1)] $\beta(\fs)$ is a wide subcategory of $\modd\Lambda$.
		\item[(2)] $\fs$ is widely generated if and only if $\fs=\tilde{\F}(\beta(\fs))$.
		\item[(3)] For any wide subcategory $\mathcal{W}\subseteq\modd\Lambda$, we have $\mathcal{W}=\beta(\tilde{\F}(\mathcal{W}))$
		\item[(4)] If $C$ is the associated cosilting module, $\mathcal{S}_C\cap\modd\Lambda$ is the set of simple objects in the abelian length category $\beta(\fs)$.
	\end{itemize}
\begin{proof}
	$(1)$ and $(3)$ was proved in \cite{MS} and $(2)$ was proved in \cite{AP}.
	
	Because by \cite[Proposition 3.11]{S}, $\mathcal{S}_C\cap\modd\Lambda$ is the set of torsion free, almost torsion modules for $(\ts,\fs)$, the proof of \cite[Proposition 3.10]{AS} carries over to prove $(4)$.
\end{proof}
\end{proposition}

\begin{lemma}\label{l3.8}
	Let $\mathcal{S}$ be a semibrick in $\modd\Lambda$. Then $\mathcal{S}$ is the set of all torsion free, almost torsion modules for $\tilde{\F}(\mathcal{S})$.
	\begin{proof}
		$\mathcal{S}$ is the set of simple objects in wide subcategory $\mathcal{W}=filt(\mathcal{S})\subseteq\modd\Lambda$. By Proposition \ref{p3.7}$(3)$, $\mathcal{W}=\beta(\tilde{\F}(\mathcal{W}))=\beta(\tilde{\F}(\mathcal{S}))$. Thus the result follow from Proposition \ref{p3.7}$(4)$.
	\end{proof}
\end{lemma}

\begin{lemma}\label{l3.9}
	Let $C\in \Cosilt_w \Lambda$. Then $\fs_C=\tilde{\F}(\mathcal{S}_C\cap\modd\Lambda)$.
	\begin{proof}
		By Proposition \ref{p3.7}$(2)$, we have that $\fs_C=\tilde{\F}(\beta(\fs_C))$. On the other hand because $\beta(\fs_C)$ is an abelian length category, any object in $\beta(\fs_C)$ is a finite filtration of simple objects. Thus the result follow from Proposition \ref{p3.7}$(4)$.
	\end{proof}
\end{lemma}
	
	\begin{theorem}\label{t3.10}
		The map $\alpha$ restricts to a bijection
		\begin{align*}
		\alpha_2:\Cosilt_w \Lambda&\longrightarrow \sbrick \Lambda\\
		C&\longmapsto \mathcal{S}_C\cap \modd\Lambda
		\end{align*}
		between cosilting module $C$ such that $\fs_C$ is widely generated and the set of semibricks in $\modd\Lambda$.
		\begin{proof}
			Let $C$ be a cosilting module. By Theorem \ref{t2.24} and \cite[Proposition 3.11]{S}, $\mathcal{S}_C\cap \modd\Lambda$ is the set of all finitely presented torsion free, almost torsion modules for\\ $(^{\bot_0}C, \Cogen(C))$. So we have a well defined map $\alpha_2:\Cosilt_w \Lambda\longrightarrow \sbrick \Lambda$.
			
			First we show that it is an injective map. Let $C,C'\in \Cosilt_w\Lambda$ and $\mathcal{S}_C\cap\modd\Lambda=\mathcal{S}_{C'}\cap\modd\Lambda$. Then by Lemma \ref{l3.9}, $\fs_C=\fs_{C'}$. Now by Remark \ref{r2.10}, $C$ and $C'$ are equivalent.
			
			Next we show that $\alpha_2$  is surjective. Let $\mathcal{S}$ be a semibrick in $\modd\Lambda$. By Theorem \ref{t2.9}, $\underrightarrow{\Lim}\tilde{\F}(\mathcal{S})=\Cogen(C)$ for some cosilting module $C$. By Lemma \ref{l3.8}, $\mathcal{S}=\mathcal{S}_C\cap\modd\Lambda$.
			
		\end{proof}
	\end{theorem}
	
	In the rest of the paper we will show that the map $\alpha$ restricts to the Asai's correspondence.
	
	\begin{lemma}\label{l3.11}
		A torsion free class $\fs\subseteq \modd\Lambda$ is functorially finite if and only if the associated cosilting module is equivalent to a support $\tau^-$-tilting module.
		\begin{proof}
			Let $M\in \modd\Lambda$ be a support $\tau^-$-tilting module. Then we can easily see that $\Cogen(M)\cap\modd\Lambda=\cogen(M)$ and it is a functorially finite torsion free class by Theorem \ref{t2.5}.
			
			Now let $C$ be a cosilting module such that $\Cogen(C)\cap\modd\Lambda\subseteq \modd\Lambda$ is functorially finite. Then by Theorem \ref{t2.5}, there is a support $\tau^-$-tilting module $M$ such that $\Cogen(C)\cap\modd\Lambda=\cogen(M)=\Cogen(M)\cap\modd\Lambda$. Then by Remark \ref{r2.10}, $C$ and $M$ are equivalent.
		\end{proof}
	\end{lemma}

	\begin{lemma}\label{l3.12}
		Let $\fs\subseteq \modd\Lambda$ be a functorially finite torsion free class with the associated support $\tau^-$-tilting module $M$. Then $\fs=\tilde{\F}(\mathcal{S}_M)$.
		\begin{proof}
			By the dual of \cite[Proposition 3.9]{MS} we have that $\fs=\tilde{\F}(\beta(\fs))$. Any object in $\beta(\fs)$, as a length abelian category, is a finite filtration of simple objects. By Proposition \ref{p3.7}, simple objects of $\beta(\fs)$ are $\mathcal{S}_M\cap\modd\Lambda=\mathcal{S}_M$. Thus the result follows.
		\end{proof}
	\end{lemma}

	\begin{theorem}\label{t3.13}
		The map $\alpha$ restricts to a bijection
		\begin{align*}
		\alpha_3:s\tau^--\tilt &\Lambda\longrightarrow \sbrick_r \Lambda\\
		&M\longmapsto \mathcal{S}_M 
		\end{align*}
		between support $\tau^-$-tilting modules and the set of right finite semibricks in $\modd\Lambda$.
		\begin{proof}
			Let $M\in \modd\Lambda$ be a support $\tau^-$-tilting module. By Lemma \ref{l3.12}, $\fs_M=\tilde{\F}(\mathcal{S}_M)$ and hence $\mathcal{S}_M$ is right finite. Thus $\alpha$ restricts to an injective map $\alpha_3:s\tau^--\tilt \Lambda\rightarrow \sbrick_r \Lambda$.
			
			Next we show that this map is surjective. Let $\mathcal{S}$ be a right finite semibrick in $\modd\Lambda$. Thus $\tilde{\F}(\mathcal{S})\subseteq\modd\Lambda$ is a functorially finite torsion free class. Therefor there exists a support $\tau^-$-tilting module $M$ such that $\tilde{\F}(\mathcal{S})=\cogen(M)=\Cogen(M)\cap\modd\Lambda$. Then by Lemma \ref{l3.8} and \cite[Proposition 3.11]{S} $\mathcal{S}_M=\mathcal{S}$.
		\end{proof}
	\end{theorem}

	\section*{acknowledgements}
	The first author would like to thanks Lidia Angeleri H{\"u}ugel and Francesco Sentieri for helpful discussions. The research of the first author was supported by a grant from IPM. The research of the second author was in part supported by a grant from IPM (No. 1403160416). The work of the second author is
based upon research funded by Iran National Science Foundation (INSF) under project No. 4001480.

\end{document}